\documentclass[11pt]{article}

\usepackage{latexsym}
\usepackage{amssymb}
\usepackage{graphicx}

\newtheorem{Theorem}{Theorem}[part]
\newtheorem{Definition}{Definition}[part]

\newtheorem{Lemma}{Lemma}[part]

\newtheorem{Remark}{Remark}[part]

\renewcommand{\theLemma}{\thesection.\arabic{Lemma}}

\renewcommand{\theequation}{\thesection.\arabic{equation}}

\def \I{\mathbb{I}}
\def \N{\mathbb{N}}
\def \R{\mathbb{R}}

\def \F{\mathbb{F}}

\def \Ac{{\cal A}}

\def \Cc{{\cal C}}

\def \Fc{{\cal F}}

\def \Lc{{\cal L}}

\def \Oc{{\cal O}}

\def \Sc{{\cal S}}

\def \Xc{{\cal X}}

\def \eps{\varepsilon}

\def \ep{\hbox{ }\hfill$\Box$}

\def\reff#1{{\rm(\ref{#1})}}

\def\beqs{\begin{eqnarray*}}
\def\enqs{\end{eqnarray*}}
\def\beq{\begin{eqnarray}}
\def\enq{\end{eqnarray}}

\def\vsp#1{\vspace{#1}}

\begin{document}

\title{On the smooth-fit property for one-dimensional \\ 
optimal switching problem}

\author{Huy{\^e}n PHAM
             \\\small  Laboratoire de Probabilit{\'e}s et
             \\\small  Mod{\`e}les Al{\'e}atoires
             \\\small  CNRS, UMR 7599
             \\\small  Universit{\'e} Paris 7
             \\\small  e-mail: pham@math.jussieu.fr
             \\\small  and CREST}

\maketitle

\begin{abstract}
This paper studies the problem of optimal switching for one-dimensional diffusion, which may be regarded as 
sequential optimal stopping problem with changes of regimes.  The resulting dynamic 
programming principle leads to a system of variational inequa\-lities, and the state space 
is divided into continuation regions and switching regions. By means of viscosity solutions 
approach, we prove the smoot-fit $C^1$ property of the value functions. 
\end{abstract}

\vspace{7mm}

\noindent {\bf Key words~:} Optimal switching, system of variational inequalities, viscosity solutions, 
smooth-fit principle.

\vspace{5mm}

\noindent {\bf MSC Classification (2000)~:} 60G40, 49L25, 60H30.

\newpage

\section{Introduction}

\setcounter{equation}{0}
\setcounter{Assumption}{0}
\setcounter{Theorem}{0}
\setcounter{Proposition}{0}
\setcounter{Corollary}{0}
\setcounter{Lemma}{0}
\setcounter{Definition}{0}
\setcounter{Remark}{0}

In this paper, we consider the optimal switching problem for a one dimensional stochastic process $X$. 
The diffusion process $X$ may take a finite number of regimes that are switched at time decisions. 
The evolution of the controlled system is  governed by 
\beqs
dX_t &=& b(X_t,I_t) dt +  \sigma(X_t,I_t) dW_t,
\enqs
with the indicator process of the regimes~:
\beqs
I_t &=& \sum_n \kappa_n 1_{\tau_n \leq t < \tau_{n+1}}. 
\enqs
Here $W$ is a standard Bronian motion on a filtred probability 
space $(\Omega,\Fc,\F=(\Fc_t)_{t\geq 0},P)$, $b$, $\sigma$ are given maps, 
$(\tau_n)_n$ is a sequence of  increasing stopping times  representing the switching regimes time decisions, and  
$\kappa_n$ is $\Fc_{\tau_n}$-measurable valued in a finite set, 
representing the new chosen value of the regime at time $\tau_n$ and until $\tau_{n+1}$. 

Our problem consists in maximizing over the switching controls $(\tau_n,\kappa_n)$ the gain functional
\beqs
E\left[ \int_0^\infty e^{-\rho t} f(X_t,I_t) dt -
\sum_n e^{-\rho\tau_n} g_{\kappa_{n-1},\kappa_n}
\right]
\enqs
where $f$ is some  running profit  function depending on the current state and the regime, 
and $g_{ij}$ is the cost for switching from regime $i$ to $j$. We then 
denote by $v_i(x)$ the value function for this control problem when starting initially from state $x$ and regime $i$.

Optimal switching problems for stochastic systems were studied by several authors, see \cite{benlio82}, \cite{lenbel83} 
or \cite{tanyon93}. These control problems lead via the dynamic programming principle to 
a system of second order variational inequalities  for the value functions $v_i$. Since the $v_i$ are not smooth $C^2$ in 
general, a first mathematical point is to give a rigorous meaning to these variational PDE, 
either in Sobolev spaces as in 
\cite{lenbel83}, or by means of viscosity notion as in \cite{tanyon93}.  
We also see that for each fixed regime $i$, the state space is divided into a switching region where it 
is optimal to change from regime $i$ to some regime $j$, and the continuation region 
where it is optimal to stay in the current regime $i$.  Optimal switching problem may 
be viewed as sequential optimal stopping problems with regimes shifts. It is well-known that optimal stopping 
problem leads to a free-boundary problem related to a variational inequality that divides the state space into 
the stopping region and the continuation region. Moreover, there is the so-called 
smooth-fit principle for optimal stopping 
problems that states the smoothness $C^1$ regularity of the value function through the boundary of the stopping region, 
once the reward function is smooth $C^1$ or is convex, see 
e.g. \cite{shi78}. Smooth-fit principle for optimal stopping problems may be proved by different arguments and we mention 
recent ones  in \cite{jac93} or \cite{pes02}  based on local time and extended It\^o's formula. 
Our main concern is to study such smooth-fit principle in the context of optimal switching problem, which has not yet 
been considered in the literature to the best of our knowledge.

Here, we use viscosity solutions arguments to prove the smooth-fit $C^1$ property of the value functions through the 
boundaries of the switching regions. The main difficulty with regard to optimal stopping problems, comes from the 
fact that the switching region for the value function $v_i$ depend also on the other value functions 
$v_j$ for which one does not know a priori $C^1$ regularity  (this is what we want to prove!) or convexity 
property. For this reason, it is an open question to see how extended It\^o's formula and local time may be used to derive such smooth-fit property for optimal switching problems.  Our proof arguments are relatively simple and 
does not require any specific knowledge on viscosity solutions theory.

The plan of this paper is organized as follows. In Section 2, we formulate our optimal switching problem and 
make some assumptions. Section 3 is devoted  to the dynamic programming PDE characterization 
of the value functions  by  viscosity solutions, through a system of variational inequalities. In Section 4, 
we prove the  smooth-fit property of the value functions.

\section{Problem formulation  and assumptions}

\setcounter{equation}{0}
\setcounter{Assumption}{0}
\setcounter{Theorem}{0}
\setcounter{Proposition}{0}
\setcounter{Corollary}{0}
\setcounter{Lemma}{0}
\setcounter{Definition}{0}
\setcounter{Remark}{0}

We start with the mathematical framework for our optimal switching problem. 
The stochastic system $X$  is valued in the state space 
$\Xc$ $\subset$ $\R$ assumed to be an interval with  endpoints 
$-\infty\leq \ell < r\leq \infty$. We let $\I_d$ $=$ $\{1,\ldots,d\}$ the finite set of regimes. 
The dynamics of the controlled stochastic system is modeled as follows. We are given maps $b$, $\sigma$ $:$ $\Xc\times\I_d$ $\rightarrow$ $\R$  satisfying a Lipschitz condition in $x$~:

\vsp{2mm}
 
{\bf (H1)} \hspace{9mm}  $|b(x,i)-b(y,i)| + |\sigma(x,i)-\sigma(y,i)|$  $\leq$  $C |x-y|$, 
\hspace{3mm}  $\forall x,y \in \Xc, \; i \in \I_d$, 

\vsp{2mm}

\noindent for some positive constant $C$, and we require 

\vsp{2mm}

{\bf (H2)} \hspace{9mm}   $\sigma(x,i)$   $>$   $0$, \hspace{3mm}  $\forall x \in {\rm int}(\Xc) \;= \; 
(\ell,r), 
\; i \in \I_d$.  

\vsp{2mm}

\noindent We set $b_i(.)$ $=$ $b(.,i)$, $\sigma_i(.)$ $=$ $\sigma(.,i)$, $i$ $\in$ $\I_d$, and we assume that  for any $x$ $\in$ $\Xc$, $i$ $\in$ $\I_d$, there exists a unique strong solution valued in $\Xc$ to  the s.d.e. 
\beq \label{diffX}
dX_t &=& b_i(X_t)dt+\sigma_i(X_t)dW_t, \;\;\; X_0 = x. 
\enq
where $W$ is a standard Brownian motion on  a filtered probability space
$(\Omega,\Fc,\F=(\Fc_t)_{t\geq 0},P)$ satisfying the usual conditions.

A switching control $\alpha$ consists of a double sequence
$\tau_1,\ldots,\tau_n,\ldots,\kappa_1,\ldots,\kappa_n,\ldots$, $n$
$\in$ $\N^*$ $=$ $\N\setminus\{0\}$, where $\tau_n$ are stopping times, $\tau_n$ $<$
$\tau_{n+1}$ and $\tau_n$ $\rightarrow$ $\infty$ a.s., and $\kappa_n$ is $\Fc_{\tau_n}$-measurable  valued in $\I_d$. We denote by $\Ac$ the set of all such switching 
controls.  Now, for any initial condition $(x,i)$ $\in$ $\Xc\times\I_d$, and any control 
$\alpha$ $=$ $(\tau_n,\kappa_n)_{n\geq 1}$ $\in$ $\Ac$, there  exists a unique strong solution valued in $\Xc\times\I_d$ to the controlled stochastic system~:
\beq
X_0 &=& x, \;\;\; I_{0^-} \; = \; i, \label{eqinitial} \\
dX_t &=& b_{_{\kappa_n}} (X_t) dt + \sigma_{_{\kappa_n}} (X_t)  dW_t,
\;\;\; I_t \; = \; \kappa_n, \;\;\;  \label{eqXn} \tau_n \leq t  <
\tau_{n+1}, \;\;\; n \geq 0. \label{eqXalpha}
\enq
Here, we set $\tau_0$ $=$ $0$ and $\kappa_0$ $=$ $i$. We denote by $(X^{x,i},I^i)$ this solution 
(as usual, we omit the dependance in $\alpha$ for notational simplicity). We notice that 
$X^{x,i}$ is a continuous process and $I^i$ is a cadlag process, possibly with a jump at time $0$ if $\tau_1$ $=$ $0$ and so $I_0$ $=$ $\kappa_1$. 

We are given a running profit function $f$ $:$ $\Xc\times\I_d$ $\rightarrow$ $\R$, and we assume  a Lipschitz condition~: 

\vsp{2mm}

{\bf (H3)} \hspace{9mm} $|f(x,i)-f(y,i)|$ $\leq$  $C |x-y|$, \hspace{3mm}  $\forall x,y \in \Xc, \; i \in \I_d$, 

\vsp{2mm}

\noindent for some positive constant $C$. We also set $f_i(.)$ $=$ $f(.,i)$, $i$ $\in$ $\I_d$.  
The cost for switching from regime $i$ to  $j$  is constant equal to  
$g_{ij}$. We assume that~: 

\vsp{2mm}

{\bf (H4)} \hspace{9mm}
$0$ $<$ $g_{ik}$  $\leq$ $g_{ij}  + g_{jk}$, \hspace{3mm} 
$\forall  i \neq j \neq k \neq i \; \in \I_d$.

\vsp{2mm}

\noindent This last condition means that the switching cost is positive and it is no more expensive to 
switch directly in one step from regime $i$ to $k$ than in two steps via an intermediate regime $j$.

The expected total profit of running the system when initial
state is $(x,i)$ and using the switching control $\alpha$ $=$
$(\tau_n,\kappa_n)_{n\geq 1}$ $\in$ $\Ac$ is
\beqs J(x,i,\alpha) &=&
E\left[ \int_0^\infty e^{-\rho t} f(X_t^{x,i},I_t^i) dt -
\sum_{n=1}^\infty e^{-\rho\tau_n} g_{\kappa_{n-1},\kappa_n} 
\right],
\enqs
where $\kappa_0$ $=$ $i$. Here $\rho$ $>$ $0$ is a positive discount factor,
and we use the convention that  $e^{-\rho\tau_n(\omega)}$  $=$ $0$
when $\tau_n(\omega)$ $=$ $\infty$. We shall see below in Lemma \ref{lemvigrowth} that 
the expectation defining $J(x,i,\alpha)$ is well-defined for $\rho$ large enough (independent of 
$x,i,\alpha$). 
The objective is to maximize this expected total profit over all
strategies $\alpha$. Accordingly, we define the function
\beq \label{defv}  v(x,i) &=& \sup_{\alpha\in\Ac} J(x,i,\alpha),
\;\;\; x \in \Xc, \; i \in \I_d.
\enq
and we denote $v_i(.)$ $:=$ $v(.,i)$ for $i$ $\in$ $\I_d$. The goal of this paper is to study the smoothness 
property of the value functions $v_i$. Our main result is the following~:

\begin{Theorem}
Assume that {\bf (H1)}, {\bf (H2)}, {\bf (H3)} and {\bf (H4)} hold.  
Then, for all $i$ $\in$ $\I_d$, the value function $v_i$ is continuously differentiable on ${\rm int}(\Xc)$ $=$ $(\ell,r)$.
\end{Theorem}

\section{Dynamic programming, viscosity solutions and system of variational inequa\-lities}

\setcounter{equation}{0}
\setcounter{Assumption}{0}
\setcounter{Theorem}{0}
\setcounter{Proposition}{0}
\setcounter{Corollary}{0}
\setcounter{Lemma}{0}
\setcounter{Definition}{0}
\setcounter{Remark}{0}

We first show the  Lipschitz continuity of the value functions $v_i$.

\begin{Lemma} \label{lemvigrowth}
Under {\bf (H1)} and {\bf (H3)}, 
there exists some positive constant $C$ $>$ $0$ such that for all $\rho$ $\geq$ $C$, 
we have~: 
\beq \label{inegvi}  
|v_i(x)-v_i(y)| & \leq & C |x-y|, \;\;\; \forall x,y \in \Xc, \; i\in \I_d.  
\enq 
\end{Lemma}
{\bf Proof.}  In the sequel, for notational simplicity, the $C$ denotes a generic constant in different places, 
depending on the constants appearing in the Lipschitz conditions in {\bf (H1)} and 
{\bf (H3)}. For any $\alpha$ $\in$ $\Ac$, the solution to \reff{eqinitial}-\reff{eqXalpha} is written as~:
\beqs
X_t^{x,i} &=& x + \int_0^t b(X_s^{x,i},I_s^i) ds + \int_0^t \sigma(X_s^{x,i},I_s^i) dW_s \\
I_t^i &=& \sum_{n=0}^\infty \kappa_n 1_{\tau_n\leq t < \tau_{n+1}}, \;\;\; (\tau_0 = 0, \; \kappa_0 =i). 
\enqs
By standard estimate for s.d.e. applying It\^o's formula to $|X_t^{x,i}|^2$ and using Gronwall's lemma, we then obtain 
from the linear growth condition on $b$ and $\sigma$ in {\bf (H1)} 
the following inequality for any $\alpha$ $\in$ $\Ac$~: 
\beqs
E \left| X_t^{x,i} \right|^2 & \leq & C e^{Ct}(1+ |x|^2), \;\;\; t \geq 0. 
\enqs
Hence, by  linear growth condition on $f$ in {\bf (H3)}, this proves that for any $\alpha$ $\in$ $\Ac$~:
\beq
E \left[ \int_0^\infty e^{-\rho t} \left|f(X_t^{x,i},I_t^i)\right| dt \right] & \leq & 
C E \left[ \int_0^\infty e^{-\rho t} (1 + |X_t^{x,i}|) dt \right] \nonumber \\
& \leq & C \int_0^\infty e^{-\rho t} e^{Ct} (1+|x|) dt \nonumber \\
&\leq & C(1+|x|),   \nonumber
\enq
for $\rho$ larger than $C$. Recalling that the $g_{ij}$ are nonnegative, 
this last inequality proves in particular that for all $(x,i,\alpha)$ $\in$ $\Xc\times\I_d\times\Ac$, 
$J(x,i,\alpha)$ is well-defined, valued in $[-\infty,\infty)$.

Moreover, by  standard estimate for s.d.e. applying It\^o's formula to $|X_t^{x,i}-X_t^{y,i}|^2$ and using 
Gronwall's lemma, we then obtain from the Lipschitz  condition {\bf (H1)} 
the following inequality uniformly in  $\alpha$ $\in$ $\Ac$~: 
\beqs
E \left| X_t^{x,i} - X_t^{y,i} \right|^2 & \leq & e^{Ct}|x-y|^2, \;\;\; \forall x,y \in\Xc, \; t \geq 0. 
\enqs
From the Lipschitz condition {\bf (H3)}, we deduce
\beqs
|v_i(x) - v_i(y)| &\leq & \sup_{\alpha\in\Ac} 
E \left[ \int_0^\infty e^{-\rho t} \left|f(X_t^{x,i},I_t^i) - f(X_t^{y,i},I_t^i) \right| dt \right] \\
&\leq & C \sup_{\alpha\in\Ac} E \left[ \int_0^\infty e^{-\rho t} \left| X_t^{x,i} - X_t^{y,i}\right| dt \right] \\
&\leq & C \int_0^\infty e^{-\rho t} e^{Ct} |x-y| dt \; \leq \; C |x-y|, 
\enqs
for $\rho$ larger than $C$.  This proves \reff{inegvi}. 
\ep

\vspace{3mm}

In the rest of this paper, we shall now assume that $\rho$ is large enough so that from the previous Lemma, 
the expected gain functional $J(x,i,\alpha)$ is  well-defined for all $x,i,\alpha$, 
and also the value functions $v_i$ are continuous.

\vspace{3mm}

The dynamic programming principle is a well-known property in stochastic optimal control. 
In our optimal switching control problem, it is formulated as follows~:

\vspace{1mm}

\noindent {\sc Dynamic programming principle~:} For any $(x,i)$ $\in$ $\Xc\times\I_d$, we have
\beq \label{progdyn}
v(x,i) &=& \sup_{_{\tiny{(\tau_n,\kappa_n)_{_n} \in\Ac}}} E\left[ \int_0^\theta  e^{-\rho t} f(X_t^{x,i},I_t^i) dt + e^{-\rho\theta} 
v(X_\theta^{x,i},I_\theta^{i}) -
\sum_{\tau_n \leq \theta}  e^{-\rho\tau_n} g_{\kappa_{n-1},\kappa_n} 
\right],
\enq 
where $\theta$ is  any stopping time, possibly depending on $\alpha$ $\in$ $\Ac$ in \reff{progdyn}. 
This principle was formally stated in \cite{benlio82} and proved rigorously for the finite horizon case in 
\cite{tanyon93}. 
The arguments for the infinite horizon case may be adapted in a straightforward way.

\vsp{3mm}

The dynamic programming principle combined with the notion of viscosity solutions are known to be a general and 
powerful tool for characterizing the value function of a stochastic control problem via a PDE representation, see  
\cite{fleson93}.  We recall the definition of viscosity solutions for a P.D.E in the form 
\begin{eqnarray} \label{FEDP}
F(x,v,D_x v, D_{xx}^2 v) &=& 0, \;\;\;  x \in \;   \Oc, 
\end{eqnarray}
where $\Oc$ is an open subset in $\R^n$ and $F$ is a continuous function and noninceasing 
in its last argument (with respect to the order of symmetric matrices).

\begin{Definition}
Let $v$ be a continuous function on $\Oc$. We say that $v$ is a viscosity solution to \reff{FEDP} on $\Oc$ it it is 

\vspace{1mm}

\noindent (i) a viscosity supersolution to \reff{FEDP} on $\Oc$~: for any $x_0$ $\in$ $\Oc$ and any $C^2$ 
function $\varphi$ in a neighborhood of $x_0$ s.t. $x_0$ is a local minimum of $v-\varphi$ and
$(v-\varphi)(x_0)$ $=$ $0$, we have~:
\begin{eqnarray*}
F(x_0,\varphi(x_0),D_x \varphi(x_0), D_{xx}^2 \varphi(x_0)) &\geq & 0. 
\end{eqnarray*}
and 

\vspace{1mm}

\noindent (ii) a viscosity subsolution to \reff{FEDP} on $\Oc$~: for any $x_0$ $\in$ $\Oc$ and any $C^2$ 
function $\varphi$ in a neighborhood of $x_0$ s.t. $x_0$ is a local maximum of $v-\varphi$ and
$(v-\varphi)(x_0)$ $=$ $0$, we have~:
\begin{eqnarray*}
F(x_0,\varphi(x_0),D_x \varphi(x_0), D_{xx}^2 \varphi(x_0)) & \leq & 0. 
\end{eqnarray*}
\end{Definition}

We shall denote by $\Lc_i$ the second order operator on the interior $(\ell,r)$ of 
$\Xc$ associated to the diffusion $X$ solution to \reff{diffX}~:  
\beqs 
\Lc_i \varphi &=& \frac{1}{2}
\sigma_i^2  \varphi" + b_i   \varphi', \;\;\; i \in \I_d. 
\enqs

\begin{Theorem} \label{thmvisco}
Assume that {\bf (H1)} and {\bf (H3)} hold. Then, for each $i$ $\in$ $\I_d$, the value function $v_i$ 
is a continuous viscosity solution on $(\ell,r)$ to the variational inequality~:
\beq \label{HJBsys}
\min\left\{ \rho v_i - \Lc_i v_i - f_i \; , \; v_i - \max_{j\neq i} (v_j - g_{ij}) \right\} &=& 0, 
\;\;\; x \in (\ell,r). 
\enq
This means that for all $i$ $\in$ $\I_d$, we have both supersolution and subsolution properties~:

\noindent (1) {\it Viscosity supersolution property~:} for any $\bar x$ $\in$ $(\ell,r)$ and $\varphi$ $\in$ $C^2(\ell,r)$ s.t. 
$\bar x$ is a local minimum  of $v_i-\varphi$, $v_i(\bar x)$ $=$ $\varphi(\bar x)$, we have 
\beq \label{supersol}
\min\left\{ \rho \varphi(\bar x) - \Lc_i \varphi(\bar x) - f_i(\bar x) \; , 
\; v_i(\bar x) - \max_{j\neq i} (v_j - g_{ij})(\bar x) \right\} & \geq & 0,
\enq 

\noindent (2) {\it Viscosity subsolution property~:} for any $\bar x$ $\in$ $(\ell,r)$ and $\varphi$ $\in$ $C^2(\ell,r)$ s.t. 
$\bar x$ is a local maximum  of $v_i-\varphi$, $v_i(\bar x)$ $=$ $\varphi(\bar x)$, we have 
\beq \label{subsol}
\min\left\{ \rho \varphi(\bar x) - \Lc_i \varphi(\bar x) - f_i(\bar x) \; , 
\; v_i(\bar x) - \max_{j\neq i} (v_j - g_{ij})(\bar x) \right\} & \leq & 0,
\enq

\end{Theorem}
{\bf Proof.}
The arguments of this proof are standard, based on the dynamic programming principle and It\^o's formula. 
We defer the proof in the appendix.  
\ep

\vspace{3mm}

For any regime $i$ $\in$ $\I_d$, we introduce the switching region~:
\beqs
\Sc_i &=& \left\{ x \in (\ell,r)~: v_i(x) = \max_{j\neq i} (v_j - g_{ij})(x) \right\}.
\enqs
$\Sc_i$ is a closed subset of $(\ell,r)$ and corresponds to the region where it is
optimal to change of regime. The complement set $\Cc_i$ of $\Sc_i$ in
$(\ell,r)$ is the so-called  continuation region~:
\beqs
\Cc_i &=& \left\{ x \in (\ell,r)~: v_i(x)  >  \max_{j\neq i} (v_j - g_{ij})(x) \right\},
\enqs
where one remains in regime $i$.

\vspace{3mm}

\begin{Remark}
{\rm Let us consider the following optimal stopping problem~: 
\beq \label{optsto}
v(x) &=& \sup_{\tau \mbox{ \tiny{stopping times}} } E\left[ \int_0^\tau e^{-\rho\tau} f(X_t^x)dt + e^{-\rho\tau} h(X_\tau^x) \right]. 
\enq
It is well-know that the dynamic programming principle for \reff{optsto} leads to a variational inequality 
for $v$ in the form~: 
\beqs
\min \left\{ \rho v - \Lc v - f \; , \; v - h \right\} &= & 0, 
\enqs
where $\Lc$ is the infinitesimal generator of the diffusion $X$. Moreover, 
the state space domain of $X$ is divided into the stopping region 
\beqs
\Sc &=& \left\{ x~: v(x) = h(x) \right\},
\enqs
and its complement set, the continuation region~: 
\beqs
\Cc &=& \left\{ x~: v(x) > h(x) \right\}.
\enqs
The smooth-fit principle for optimal stopping problems states that 
the value function $v$ is smooth $C^1$ through the boundary of the stopping region, 
the so-called free boundary, once $h$ is $C^1$ or convex. 

Our aim is to state similar results for optimal switching problems. 
The main difficulty comes from the fact that we have a system a variational 
inequalities, so that the switching region for $v_i$ depend also on the other value functions $v_j$ which are not convex or known to be $C^1$ a priori.    
}
\end{Remark}

\section{The smooth-fit property}

\setcounter{equation}{0}
\setcounter{Assumption}{0}
\setcounter{Theorem}{0}
\setcounter{Proposition}{0}
\setcounter{Corollary}{0}
\setcounter{Lemma}{0}
\setcounter{Definition}{0}
\setcounter{Remark}{0}

We first show, like for optimal stopping problems, that the value functions are  smooth $C^2$ in their 
continuation regions. We provide here a quick proof based on viscosity solutions arguments.

\begin{Lemma}  \label{lemC2Cci}
Assume that {\bf (H1)}, {\bf (H2)} and {\bf (H3)} hold. 
Then, for all $i$ $\in$ $\I_d$, the value function $v_i$ is smooth $C^2$ on $\Cc_i$ and satisfies in a classical sense~:  
\beq \label{edpcci}
\rho v_i (x) - \Lc_i v_i (x)  - f_i (x)  &=& 0, \;\;\; x \in \Cc_i.
\enq
\end{Lemma}
{\bf Proof.} We first check that $v_i$ is a viscosity solution to \reff{edpcci}. Let $\bar x$ $\in$ $\Cc_i$ and 
$\varphi$ a $C^2$ function on $\Cc_i$ s.t. $\bar x$ is a local maximum of $v_i-\varphi$, 
$v_i(\bar x)$ $=$ $\varphi(\bar x)$.  
Then, by definition of $\Cc_i$, we have $v_i(\bar x)$ $>$ $\max_{j\neq i}(v_j - g_{ij})(\bar x)$, and so from the 
subsolution viscosity property \reff{subsol} of $v_i$, we have~: 
\beqs
\rho \varphi(\bar x) - \Lc_i \varphi (\bar x)  - f_i (\bar x)  &\leq & 0. 
\enqs
The supersolution inequality for \reff{edpcci} is immediate from \reff{supersol}.

Now, for arbitrary bounded interval $(x_1,x_2)$ $\subset$ $\Cc_i$, consider the Dirichlet boundary linear problem~: 
\beq
\rho w (x) - \Lc_i w (x)  - f_i (x)  &=& 0, \;\;\; \mbox{ on } (x_1,x_2) \label{edpx12} \\
w(x_1) \; = \; v_i(x_1), & & w(x_2) \; = \; v_i(x_2). \label{edptermix12}
\enq
Under the nondegeneracy condition {\bf (H2)}, classical results provide the existence and uniqueness of a 
smooth $C^2$ function $w$ solution on $(x_1,x_2)$ to \reff{edpx12}-\reff{edptermix12}. 
In particular, this smooth function $w$ is a viscosity solution of \reff{edpcci} on $(x_1,x_2)$. From standard uniqueness 
results on viscosity solutions (here for a linear PDE in a bounded domain), we deduce that $v_i$ $=$ $w$ on $(x_1,x_2)$. 
From the arbitrariness of $(x_1,x_2)$ $\subset$ $\Cc_i$, this proves that $v_i$ is smooth $C^2$ on $\Cc_i$, and so satisfies 
\reff{edpcci} in a classical sense. 
\ep

\vspace{3mm}

We now state an elementary partition property on the switching regions.

\begin{Lemma}  \label{lemSi} 
Assume that {\bf (H1)}, {\bf (H3)} and {\bf (H4)} hold. Then, for all $i$ $\in$ $\I_d$, we have $\Sc_i$ $=$ $\cup_{j\neq i} \Sc_{ij}$ where
\beqs
\Sc_{ij} &=& \left\{ x \in \Cc_j~: v_i(x) = (v_j - g_{ij})(x) \right\}.
\enqs
\end{Lemma}
 {\bf Proof.}  Denote $\tilde \Sc_i$ $=$ $\cup_{j\neq i} \Sc_{ij}$. Since we always have 
$v_i$ $\geq$ $\max_{j\neq i}(v_j-g_{ij})$, the inclusion $\tilde \Sc_i$ $\subset$ $\Sc_i$ is clear. 

Conversely, let $x$ $\in$ $\Sc_i$. Then there exists $j$ $\neq$ $i$ s.t. $v_i(x)$ $=$ 
$v_j(x) - g_{ij}$. We have two cases~:  

$\star$ if $x$ lies in  $\Cc_j$, then $x$ $\in$ $\Sc_{ij}$ and so $x$ $\in$ $\tilde \Sc_i$.

$\star$ if $x$ does not lie in $\Cc_j$, then $x$ would lie in $\Sc_j$, which means 
that one could find some $k$ $\neq$ $j$ s.t. $v_j(x)$ $=$ $v_k(x)-g_{jk}$, and so  
$v_i(x)$ $=$ $v_k(x)-g_{ij}-g_{jk}$. From condition {\bf (H4)}  and since we always have 
$v_i$ $\geq$ $v_k-g_{ik}$, this would imply $v_i(x)$ $=$ $v_k(x)-g_{ik}$. Since cost $g_{ik}$ is positive, 
we also get that $k$ $\neq$ $i$. Again, we have two cases~: if  $x$ lies in $\Cc_k$, then $x$ lies in 
$\Sc_{ik}$ and so in $\tilde\Sc_i$. Otherwise, we repeat the above argument and since the number of states is 
finite, we should necessarily find some  $l$ $\neq$ $i$ s.t. $v_i(x)$ $=$ $v_l(x)- g_{il}$ and 
$x$ $\in$ $\Cc_l$. This shows finally that $x$ $\in$ $\tilde\Sc_i$.  
\ep

\vspace{2mm}

\begin{Remark}
{\rm $\Sc_{ij}$ represents the region where it is optimal  to switch from regime $i$ to 
regime $j$ and stay here for a moment, i.e. without changing instantaneously from regime $j$ to another regime.
}
\end{Remark}

\vspace{2mm}

We can finally  prove the smooth-fit property of the value functions $v_i$ through the boundaries of the switching regions.

\begin{Theorem} \label{theosmoothfit} 
Assume that {\bf (H1)}, {\bf (H2)}, {\bf (H3)} and {\bf (H4)} hold.  
Then, for all $i$ $\in$ $\I_d$, the value function $v_i$ is continuously differentiable on $(\ell,r)$. Moreover, 
at $x$ $\in$ $\Sc_{ij}$, we have $v_i'(x)$ $=$ $v_j'(x)$.  
\end{Theorem}
{\bf Proof.}  We already know from Lemma \ref{lemC2Cci} that $v_i$ is smooth $C^2$ on the open set $\Cc_i$ for all $î$ $\in$ $\I_d$. We have 
to prove the $C^1$ property of $v_i$ at any point of the closed set $\Sc_i$. 
We denote for all $j$ $\in$ $\I_d$, $j$ $\neq$ $i$, $h_j$ $=$ $v_j - g_{ij}$ and we notice that $h_j$ 
is smooth $C^1$ (actually even $C^2$) on $\Cc_j$.

\vspace{2mm}

\noindent {\bf 1.} We first check that $v_i$ admits a left and right derivative 
$v'_{i,-} (x_0)$ and $v'_{i,+} (x_0)$ at any point $x_0$ in $\Sc_i$ $=$ $\cup_{j\neq i}\Sc_{ij}$. We distinguish 
the two following cases~: 

\vspace{1mm}

\noindent $\bullet$ {\it Case a)} $x_0$ lies in the interior  ${\rm Int}(\Sc_i)$ of $\Sc_i$. Then, we have two 
subcases~: 

\vspace{1mm}

$\star$ $x_0$ $\in$  ${\rm Int}(\Sc_{ij})$ for some $j$ $\neq$ $i$, i.e. there exists some $\delta$  $>$ $0$ 
s.t. $[x_0-\delta,x_0+\delta]$ $\subset$ $\Sc_{ij}$. By definition of $\Sc_{ij}$, we then have $v_i$ $=$ $h_j$ on 
$[x_0-\delta,x_0+\delta]$ $\subset$ $\Cc_{j}$, and so $v_i$ is differentiable at $x_0$ with 
$v'_i(x_0)$ $=$ $h_j'(x_0)$. 

\vspace{1mm}

$\star$ There exists $j\neq k\neq i$ in $\I_d$ and  $\delta$ $>$ $0$  s.t. $[x_0-\delta,x_0]$  
$\subset$ $\Sc_{ij}$ and $[x_0,x_0+\delta]$  $\subset$ $\Sc_{ik}$. We then have  
$v_i$ $=$ $h_j$ on $[x_0-\delta,x_0]$ $\subset$ $\Cc_{j}$ and $v_i$ $=$ $h_k$ on 
$[x_0,x_0+\delta]$ $\subset$ $\Cc_{k}$. Thus, $v_i$ admits a left and right derivative at $x_0$ with 
$v'_{i,-} (x_0)$ $=$ $h_j'(x_0)$ and $v'_{i,+} (x_0)$ $=$ $h_k'(x_0)$.  

\vspace{1mm}

\noindent $\bullet$ {\it Case b)} $x_0$ lies in the boundary $\partial\Sc_i$ $=$ 
$\Sc_i\setminus{\rm Int}(\Sc_i)$ of $\Sc_i$. We assume that $x_0$ lies in the left-boundary of $\Sc_i$, i.e. 
there exists $\delta$ $>$ $0$ s.t. $[x_0-\delta,x_0)$ $\subset$ $\Cc_i$ (the other case where $x_0$ lies in the 
right-boundary is dealt with similarly).  Recalling that on $\Cc_i$, 
$v_i$ is solution to~: $\rho v_i-\Lc v_i-f_i$ $=$ $0$, we deduce that on $[x_0-\delta,x_0)$, $v_i$ is equal 
to $w_i$  
the unique smooth $C^2$ solution to the o.d.e.~: $\rho w_i - \Lc w_i - f_i$ $=$ $0$ 
with the boundaries conditions~:  
$w_i(x_0-\delta)$ $=$ $v_i(x_0-\delta)$, $w_i(x_0)$ $=$ $v_i(x_0)$. Therefore, $v_i$ admits a left 
derivative at $x_0$ with $v'_{i,-} (x_0)$ $=$ $w_i'(x_0)$. In order to prove that $v_i$ admits a right derivative, 
we distinguish the two subcases~: 

\vspace{1mm}

$\star$  There exists $j\neq i$ in $\I_d$ and  $\delta'$ $>$ $0$  s.t. $[x_0,x_0+\delta']$ 
$\subset$ $\Sc_{ij}$. Then, on $[x_0,x_0+\delta']$, $v_i$ is 
equal to $h_j$. Hence $v_i$ admits a  right derivative at $x_0$ with $v'_{i,+} (x_0)$ $=$ $h_j'(x_0)$.

$\star$ Otherwise, for all $j$ $\neq$ $i$, we can find a sequence 
$(x_n^j)$ s.t. $x_n^j$ $\geq$ $x_0$, $x_n^j$ $\notin$ $\Sc_{ij}$ and $x_n^j$ $\rightarrow$ 
$x_0$. By a diagonalization procedure, we construct then a sequence $(x_n)$ s.t. $x^n$ $\geq$ $x_0$, 
$x_n$ $\notin$ $\Sc_{ij}$ for all $j$ $\neq$ $i$, i.e. $x_n$ $\in$ $\Cc_i$, and 
$x_n$ $\rightarrow$ $x_0$. Since $\Cc_i$ is open, there exists then $\delta''$ $>$ $0$ s.t. 
$[x_0,x_0+\delta'']$ $\subset$ $\Cc_i$. 
We deduce that on $[x_0,x_0+\delta'']$, $v_i$ is equal to $\hat w_i$ the unique smooth $C^2$ solution to the o.d.e. 
$\rho \hat w_i - \Lc\hat w_i - f_i$ $=$ $0$ with the boundaries conditions 
$\hat w_i(x_0)$ $=$ $v_i(x_0)$, $\hat w_i(x_0+\delta'')$ $=$ $v_i(x_0+\delta'')$. In particular, 
$v_i$ admits a right derivative at $x_0$ with $v'_{i,+} (x_0)$ $=$ $\hat w_i'(x_0)$.

\vspace{1mm}

\vspace{2mm}

\noindent {\bf 2.} Consider now some point in $\Sc_i$ eventually on its boundary. 
We recall again that from Lemma \ref{lemSi}, there exists some $j$ $\neq$ $i$ s.t. $x_0$ $\in$ $\Sc_{ij}$~: 
$v_i(x_0)$ $=$ $h_j(x_0)$, and  $h_j$ is smooth $C^1$on $x_0$ in $\Cc_j$. Since $v_j$ $\geq$ $h_j$, we deduce that 
\beqs
\frac{v_i(x)-v_i(x_0)}{x-x_0} &\leq & \frac{h_j(x)-h_j(x_0)}{x-x_0}, \;\;\; \forall \; x < x_0 \\
\frac{v_i(x)-v_i(x_0)}{x-x_0} &\geq & \frac{h_j(x)-h_j(x_0)}{x-x_0}, \;\;\; \forall \; x > x_0, 
\enqs
and so~:
\beqs
v'_{i,-} (x_0) \;  \leq & h'_j(x_0) & \leq 
v'_{i,+} (x_0). 
\enqs
We argue by contradiction and suppose that $v_i$ is not differentiable at $x_0$. Then, in view of the 
above inequality, one can find some $p$ $\in$ $(v'_{i,-} (x_0),v'_{i,+} (x_0))$. Consider, for $\eps$ $>$ $0$, 
the smooth $C^2$ function~: 
\beqs
\varphi_\eps(x) &=& v_i(x_0) + p(x-x_0) + \frac{1}{2\eps}(x-x_0)^2. 
\enqs
Then, we see that $v_i$ dominates locally in a neighborhood of $x_0$ 
the function $\varphi_\eps$, i.e $x_0$ is a local 
minimum of $v_i-\varphi_\eps$. From  the supersolution viscosity property of $v_i$ to the PDE \reff{HJBsys}, 
this yields~: 
\beqs
\rho \varphi_\eps(x_0)  - \Lc_i \varphi_\eps (x_0) - f_i(x_0) &\geq & 0, 
\enqs
which is written as~: 
\beqs
\rho v_i(x_0) - b_i (x_0) p - f_i(x_0) - \frac{1}{2\eps} \sigma_i^2 (x_0) &\geq & 0. 
\enqs
Sending $\eps$ to zero provides the required contradiction under {\bf (H2)}. We have then proved that for $x_0$ $\in$ 
$\Sc_{ij}$, $v_i'(x_0)$ $=$ $h_j'(x_0)$ $=$ $v_j'(x_0)$. 
\ep

\section*{Appendix: Proof of Theorem \ref{thmvisco}}

\renewcommand{\theLemma}{A.\arabic{Lemma}}
\renewcommand{\theequation}{A.\arabic{equation}}
\setcounter{equation}{0}
\setcounter{Lemma}{0}

(1) \underline{{\it Viscosity supersolution property}}.

\vsp{1mm}

\noindent  Fix $i$ $\in$ $\I_d$. Consider any $\bar x$ $\in$ $(\ell,r)$ and $\varphi$ $\in$ $C^2(\ell,r)$ s.t. $\bar x$ is a minimum of $v_i - \varphi$ in a neighborhood $B_\eps(\bar x)$ $=$ $(\bar x-\eps,\bar x+\eps)$ of $\bar x$, $\eps$ $>$ $0$, and  $v_i(\bar x)$ $=$ $\varphi(\bar x)$. By taking the immediate switching control $\tau_1$ $=$ $0$, 
$\kappa_1$ $=$ $j$ $\neq$ $i$, $\tau_n$ $=$ $\infty$, $n$ $\geq$ 
$2$, and $\theta$ $=$ $0$ in the relation \reff{progdyn}, we obtain 
\beq \label{vivj}
v_i(\bar x) &\geq & v_j(\bar x) - g_{ij}, \;\;\; \forall j \neq i. 
\enq 
On the other hand, by taking the no-switching control $\tau_n$ $=$ $\infty$, $n$ $\geq$ $1$, i.e. 
$I_t^i$ $=$ $i$, $t$ $\geq$ $0$, $X^{\bar x,i}$ stays in regime $i$ with diffusion coefficients $b_i$ and $\sigma_i$,  
and $\theta$ $=$ $\tau_\eps$ $\wedge$ $h$, with $h$ $>$ $0$ and $\tau_\eps$ $=$ $\inf\{t\geq 0~:$ 
$X_t^{\bar x,i} \notin B_\eps(\bar x) \}$, we get from \reff{progdyn}~: 
\beqs
\varphi(\bar x) \; = \; v_i(\bar x)  & \geq & E\left[\int_0^{\theta} e^{-\rho t} f_i(X_t^{\bar x,i}) dt + 
e^{-\rho\theta} v_i (X_{\theta}^{\bar x,i}) \right ] \\
&\geq & E\left[\int_0^{\theta} e^{-\rho t} f_i(X_t^{\bar x,i}) dt + 
e^{-\rho\theta} \varphi (X_{\theta}^{\bar x,i}) \right ]
\enqs  
By applying It\^o's formula to $e^{-\rho t}\varphi(X_t^{\bar x,i})$ between $0$ and $\theta$ $=$ 
$\tau_\eps\wedge h$ and plugging into the last inequality, we obtain~:
\beqs
\frac{1}{h} E\left[\int_0^{\tau_\eps\wedge h} e^{-\rho t} 
\left( \rho\varphi - \Lc_i \varphi - f_i\right)(X_t^{\bar x,i})  \right] & \geq & 0.
\enqs
From the dominated convergence theorem, this yields by sending $h$ to zero~:
\beqs
(\rho\varphi - \Lc_i \varphi - f_i)(\bar x) &\geq & 0. 
\enqs
By combining with \reff{vivj}, we obtain the required supersolution inequality \reff{supersol}.

\vsp{3mm}

\noindent (2) \underline{{\it Viscosity subsolution property}}.

\vsp{1mm} 

\noindent Fix $i$ $\in$ $\I_d$, and consider any $\bar x$ $\in$ $(\ell,r)$ and $\varphi$ $\in$ $C^2(\ell,r)$ s.t. $\bar x$ is a maximum of $v_i - \varphi$ in a neighborhood 
$B_\eps(\bar x)$ $=$ $(\bar x-\eps,\bar x+\eps)$ of $\bar x$, $\eps$ $>$ $0$, and $v_i(\bar x)$ $=$ $\varphi(\bar x)$.  We argue by contradiction by assuming on the 
contrary that \reff{subsol} does not hold so that by continuity of $v_i$, $v_j$, $j$ $\neq$ $i$, $\varphi$ and its derivatives, there exists some  $0$ $<$ $\delta$ $\leq$ $\eps$ s.t. 
\beq
(\rho\varphi - \Lc_i \varphi - f_i)(x) & \geq &  \delta, \;\;\; \forall x \in B_\delta(\bar x) = (x-\delta,x+\delta)  \label{Lcidelta} \\
v_i(x) \; - \;  \max_{j\neq i} (v_j - g_{ij})(x) &\geq &  \delta ,  \;\;\; \forall x \in B_\delta(\bar x). 
\label{videlta}
\enq
For any $\alpha$ $=$ $(\tau_n,\kappa_n)_{n\geq 1}$ $\in$ $\Ac$, consider the exit time 
$\tau_\delta$ $=$ $\inf\{t\geq 0~:$ 
$X_t^{\bar x,i} \notin B_\delta(\bar x) \}$. By applying It\^o's formula to 
$e^{-\rho t}\varphi(X_t^{\bar x,i})$ between $0$ and $\theta$ $=$ $\tau_1\wedge \tau_\delta$, we have by noting that before $\theta$,    $X^{x,i}$ stays in regime $i$ and in the ball 
$B_\delta(\bar x)$ $\subset$ $B_\eps(\bar x)$~:  
\beq
v_i(\bar x) \; = \; \varphi(\bar x) &=& E\left[\int_0^{\theta}e^{-\rho t} (\rho\varphi - \Lc_i \varphi)(X_t^{\bar x,i}) dt + e^{-\rho\theta} \varphi(X_\theta^{\bar x,i}) \right] \nonumber \\
&\geq &  E\left[\int_0^{\theta}e^{-\rho t} (\rho\varphi - \Lc_i \varphi)(X_t^{\bar x,i}) dt + e^{-\rho\theta}  v_i(X_\theta^{\bar x,i}) \right]. \label{interitophi} 
\enq
Now, since $\theta$  $=$ $\tau_\delta\wedge \tau_1$, we have
\beqs
e^{-\rho\theta} v(X_\theta^{\bar x,i},I_\theta^i) - \sum_{\tau_n\leq\theta} g_{\kappa_{n-1}, \kappa_n}
&=& e^{-\rho\tau_1} \left(v(X_{\tau_1}^{\bar x,i},\kappa_1) - g_{i\kappa_1}\right) 
1_{\tau_1 \leq \tau_\delta}   +   e^{-\rho\tau_\delta} v_i(X_{\tau_\delta}^{\bar x,i}) 1_{\tau_\delta < \tau_1}  \\
&\leq & e^{-\rho\tau_1} \left(v_i(X_{\tau_1}^{\bar x,i}) - \delta\right) 
1_{\tau_1 \leq \tau_\delta}   +   
e^{-\rho\tau_\delta} v_i(X_{\tau_\delta}^{\bar x,i}) 1_{\tau_\delta < \tau_1} \\
&=& e^{-\rho\theta} v_i(X_{\theta}^{\bar x,i}) - \delta e^{-\rho\tau_1} 1_{\tau_1 \leq \tau_\delta},  
\enqs
where the inequality follows from \reff{videlta}. By plugging into \reff{interitophi} and using 
\reff{Lcidelta}, we get~: 
\beq
v_i(\bar x) & \geq & E\left[\int_0^{\theta}e^{-\rho t} f_i(X_t^{\bar x,i}) dt + 
e^{-\rho\theta} v(X_\theta^{\bar x,i},I_\theta^i) - \sum_{\tau_n\leq\theta} g_{\kappa_{n-1}, \kappa_n}
\right] \nonumber \\
& & \;\;\; + \; \delta \; E\left[\int_0^\theta e^{-\rho t} dt + e^{-\rho \tau_1} 1_{\tau_1 \leq \tau_\delta} \right]. 
\label{vitointer2} 
\enq 
We now claim that there exists some positive constant $c_0$ $>$ $0$ s.t.~:  
\beqs
E\left[\int_0^\theta e^{-\rho t} dt + e^{-\rho \tau_1} 1_{\tau_1 \leq \tau_\delta} \right] &\geq & c_0, \;\;\; \forall 
\alpha \in \Ac. 
\enqs
For this, we construct a smooth function $w$  s.t. 
\beq
\max\left\{ \rho w(x) - \Lc_i w (x) - 1 \; , \; w(x) - 1 \right\} & \leq &  0, \; 
\forall x \in B_\delta(\bar x) \label{lciw} \\
w(x) & = & 0, \; \forall x \in \partial B_\delta(\bar x) \; = \; \{x~: |x-\bar x| = \delta \} 
\label{w1} \\
w(\bar x) & > & 0. \label{wbarx} 
\enq
For instance, we can take the function $w(x)$ $=$ $c_0 \left(1 - \frac{|x-\bar x|^2}{\delta^2}\right)$, 
with 
\beqs
0 < c_0 & \leq & \min\left\{ \left( \rho + \frac{2}{\delta} \sup_{x\in B_\delta(\bar x)} |b_i(x)| + 
\frac{1}{\delta^2} \sup_{x\in B_\delta(\bar x)} |\sigma_i(x)|^2 \right)^{-1}\; , \; 1 \right\}. 
\enqs
Then, by applying It\^o's formula to $e^{-\rho t} w(X_t^{\bar x,i})$ between $0$ and $\theta$ $=$ 
$\tau_\delta\wedge\tau_1$, we have~:
\beqs
0 \; < \; c_0 \; = \; w(\bar x) &=& E\left[\int_0^\theta e^{-\rho t} (\rho w - \Lc_i w) (X_t^{\bar x,i}) dt + 
e^{-\rho\theta} 
w(X_\theta^{\bar x,i}) \right]  \\
& \leq & E\left[\int_0^\theta e^{-\rho t}  dt + e^{-\rho \tau_1} 1_{\tau_1 \leq \tau_\delta} \right],
\enqs 
from \reff{lciw}, \reff{w1} and \reff{wbarx}.  By plugging this last inequality (uniform  in $\alpha$)  into \reff{vitointer2}, we then obtain~:   
\beqs
v_i(\bar x) & \geq & \sup_{\alpha\in\Ac} E\left[\int_0^{\theta}e^{-\rho t} f_i(X_t^{\bar x,i}) dt + 
e^{-\rho\theta} v(X_\theta^{\bar x,i},I_\theta^i) - \sum_{\tau_n\leq\theta} g_{\kappa_{n-1}, \kappa_n}
\right]  \; + \; \delta c_0,
\enqs 
which is in contradiction with the dynamic programming principle \reff{progdyn}.


\noindent  




\vspace{9mm}


\begin{thebibliography}{}


\bibitem{benlio82} Bensoussan A. et J.L. Lions (1982)~: Contr\^ole impulsionnel et in\'equations 
variationnelles, Dunod. 





\bibitem{jac93} Jacka S. (1993)~: ``Local times, optimal stopping and semimartingales'', {\it Annals of Probability}, 
21, 329-339.

\bibitem{fleson93} Fleming W. and M. Soner (1993)~: Controlled Markov processes and viscosity solutions, Springer Verlag. 


\bibitem{lenbel83} Lenhart S. and S. Belbas (1983)~: ``A system of nonlinear partial differential equations arising in the 
optimal control of stochastic systems with switching costs'', {\it SIAM J. Appl. Math.}, 43, 465-475. 





\bibitem{pes02} Peskir G. (2002)~: ``A change of variable formula with local time on curves'', to appear in 
{\it J. Theoret. Prob.} 



\bibitem{shi78} Shiryaev A. (1978)~: Optimal stopping rules, Springer Verlag. 



\bibitem{tanyon93} Tang S. and J. Yong (1993)~: ``Finite horizon stochastic optimal switching and impulse controls 
with a viscosity solution approach'', {\it Stoch. and Stoch. Reports}, 45, 145-176. 




\end{thebibliography}
\end{document}